\documentclass[11pt,leqno]{amsart}

\usepackage[letterpaper,margin=1.08in,footskip=0.36in]{geometry}
\usepackage[T1]{fontenc}
\IfFileExists{scholax.sty}{%
  \usepackage[p,scosf,looser]{scholax}
  \usepackage[scaled=1.08,ncf,vvarbb]{newtxmath}
}{%
  \usepackage{lmodern}
  \usepackage{mathpazo}
  \usepackage{amssymb}
}
\usepackage{mathtools}
\usepackage{aliascnt}
\usepackage{microtype}
\usepackage{enumitem}
\usepackage[dvipsnames]{xcolor}
\usepackage[colorlinks,pagebackref]{hyperref}
\usepackage[nameinlink,noabbrev]{cleveref}
\usepackage{bookmark}
\usepackage{tikz-cd}
\PassOptionsToPackage{pdfusetitle,pagebackref,colorlinks}{hyperref}
\usepackage{amsfonts}
\usepackage{bookmark}
\hypersetup{citecolor=OliveGreen,linkcolor=Mahogany,urlcolor=Plum,pdfauthor={Rahul Ajit and Shikha Bhutani}}

\allowdisplaybreaks
\numberwithin{equation}{section}
\setlist[enumerate]{leftmargin=2.3em}

\newcommand{\cO}{\mathcal O}
\newcommand{\Q}{\mathbb Q}
\newcommand{\Z}{\mathbb Z}

\newcommand{\Exc}{\operatorname{Exc}}
\newcommand{\Supp}{\operatorname{Supp}}
\newcommand{\Spec}{\operatorname{Spec}}

\newcommand{\depth}{\operatorname{depth}}

\newcommand{\floor}[1]{\lfloor #1\rfloor}

\newtheorem{theorem}{Theorem}[section]
\newtheorem{proposition}[theorem]{Proposition}
\newtheorem{lemma}[theorem]{Lemma}
\newtheorem{corollary}[theorem]{Corollary}

\theoremstyle{definition}

\theoremstyle{remark}
\newtheorem{remark}[theorem]{Remark}

\crefname{mainthm}{Theorem}{Theorems}
\Crefname{mainthm}{Theorem}{Theorems}
\crefname{theorem}{Theorem}{Theorems}
\crefname{conjecture}{Conjecture}{Conjectures}
\Crefname{conjecture}{Conjecture}{Conjectures}
\Crefname{theorem}{Theorem}{Theorems}
\crefname{proposition}{Proposition}{Propositions}
\Crefname{proposition}{Proposition}{Propositions}
\crefname{lemma}{Lemma}{Lemmas}
\Crefname{lemma}{Lemma}{Lemmas}
\crefname{remark}{Remark}{Remarks}
\Crefname{remark}{Remark}{Remarks}

\title{Kov\'acs Vanishing for Excellent LC Threefolds}

\author{Rahul Ajit}
\address{Department of Mathematics\\
The University of Utah\\
Salt Lake City, UT 84112, USA}
\email{rahulajit@math.utah.edu}

\author{Shikha Bhutani}
\address{Department of Mathematics\\
Michigan State University\\
East Lansing, MI 48824, USA}
\email{bhutani4@msu.edu}

\date{\today}

\begin{document}
\begin{abstract}
We prove a vanishing theorem of S. Kov\'acs for excellent log canonical threefolds whose closed points have perfect residue fields of sufficiently large characteristics. This vanishing theorem roughly says that LC threefold pairs are not far from being rational pairs (in the sense of Koll\'ar \& Kov\'acs and Schwede \& Takagi) if one includes a correction term. We apply this vanishing to prove a depth relation of Chou and resolution-invariance of cohomologies of ``essential divisors" for isolated log canonical threefolds in large characteristic, extending works of Ishii and Fujino. Using examples constructed by Bernasconi and Cascini-Tanaka, we show that these results fail in low characteristic. Finally, we show Kov\'acs and Alexeev-Hacon's result (in char 0) that ``Irrational centers are LC centers" fails in every positive characteristic.
\end{abstract}

\maketitle
\tableofcontents
\section{Introduction}
Vanishing theorems are fundamental tools in birational algebraic geometry. Rational singularities have been the gold standard for “mild” singularities in algebraic geometry. They are extremely useful because their cohomological behavior is very similar to that of smooth varieties. There is a close relationship between MMP singularities and rational singularities. In particular, Kawamata log terminal singularities are rational and Gorenstein rational singularities are KLT.

Unfortunately, Kodaira vanishing fails in positive characteristics, and as a result, by taking cone one sees that, the GR vanishing fails when the dimension is at least 3. On a positive note, Chatzistamatiou and R\"ulling proved in \cite{CR} that GR vanishing holds for excellent, regular schemes $X$ and $Y$ in any dimension. In dimension $3$, positive and mixed characteristics $(p > 5)$, GR vanishing holds for klt singularities, proving that they are rational, see \cite{ABL, BK, hacon_witaszek_rationality}.

Log canonical singularities are a larger class of singularities than KLT and are not necessarily rational, or even Cohen-Macaulay. In \cite[Theorem 6.1 and Corollary 6.2]{Kovacs-DB-vanishing}, Kov\'acs proved a rather general vanishing theorem for Du Bois pairs, and used it to derive a vanishing theorem for log canonical pairs over complex numbers in all dimension.
More precisely, let $(X,\Delta)$ be a log canonical pair, and
$\pi\colon Y\to X$ be a log resolution, with
$K_Y+\Delta_Y=\pi^*(K_X+\Delta).$
Let $B_\pi$ be the reduced divisor consisting of the coefficient-one components of the crepant boundary. Kov\'acs vanishing asserts that
$R^i\pi_*\mathcal O_Y(-B_\pi)=0$ for all $ i>0.$
Philosophically, Kov\'acs vanishing says that log canonical singularities are not far from being rational singularities. And, their failure
to be rational is controlled by the non-klt locus. In this sense, Kov\'acs vanishing can even be thought of analogous to the definition of rational pairs with a correction term to take care of log canonicity.; compare \cite[Definition~2.80]{KollarSingularities} and
\cite{schwede_rational_2008}. Over the complex numbers, see \cite[Corollary~6.5]{Kovacs-DB-vanishing} for the $\Q$-factorial case and \cite[Theorem~1]{Fujino-Kovacsvanishing} for the general statement.

Kov\'acs vanishing is closely related to Steenbrink's characterization
of normal isolated Du Bois singularities; see
\cite[3.6]{Steenbrink},
\cite[Proposition~4.13]{DuBois}, and
\cite[Theorem~6.1]{Kovacs-Schwede-Hodge-meets-MMP}.
A weaker form of this vanishing also plays a crucial role in the
extension theorem for differential forms on log canonical spaces
proved in \cite{GKKP}; see also \cite{Graf-DifFormLC} for related
results in positive characteristic.

Our main result is the following.

\begin{theorem} (cf.\cite[Theorem 1.2]{Kovacs-DB-vanishing})\label{Thm:mainKovacs}
Let $(X,\Delta)$ be a log canonical threefold, quasi-projective over an excellent local domain admitting a dualizing complex. Assume that the closed points of $X$ have perfect residue fields. Suppose that one of the following holds.
\begin{enumerate}[label=\textnormal{(\alph*)}]
\item Every closed point has residue characteristic greater than $5$, and $(X,\Delta)$ has no zero-dimensional LC center.
\item The nonzero coefficients of $\Delta$ belong to a DCC set
$
 I\subseteq [\tfrac12,1]\cap\Q,
$
and every closed point has residue characteristic greater than $p_I$ depending only on $I$.
\end{enumerate}
Then, for every proper log resolution $\pi\colon Y\to X$, we have $R^i\pi_*\cO_Y(-B_\pi)=0$ for all $i>0$. Moreover, $R\pi_*\cO_Y(-B_\pi)$ is independent of $\pi$.
\end{theorem}

We next derive two applications of Kov\'acs vanishing. Over $\mathbb{C}$, a fundamental distinction between klt and log canonical
singularities is that klt singularities are rational, whereas log
canonical singularities are Du Bois
\cite{Kollar-Kovacs-LC-duBois-1,Kollar-Kovacs-LC-duBois-2}, but need
not be rational or even Cohen--Macaulay. It is therefore natural to ask
how far log canonical singularities are from being rational. Since
rational singularities are Cohen--Macaulay, one basic question is to
determine when a log canonical singularity is Cohen--Macaulay, or more
generally to compute its depth.

Several results are known along this direction in characteristic zero. Fujino \cite{Fujino-Kovacsvanishing} proved that if $(X,\Delta)$ is a log canonical pair of
dimension at least three and $\overline{\{x\}}$ is not a log canonical
center, then $\operatorname{depth}_x \mathcal O_X
    \geq \min\{3,\operatorname{codim}_X x\}$
, generalizing 
\cite[Lemma~3.2]{AlexeevS3}. Koll\'ar and Kov\'acs generalized this in \cite{KollarLocalKVV,KovacsIrr}; see also \cite{Hacon-Alexeev}.

In positive characteristic, however, the direct analogue of this depth
bound fails in every characteristic. 
For every prime $p>0$, \cite[Theorem~1.3]{ABP} constructed three-dimensional log canonical singularity with a closed point $x$ which
is not itself an LC center but satisfies $\operatorname{depth}_x \mathcal O_X=2.$

Suppose residue fields are perfect of characteristic greater than five. Let $x$ be a closed point on a one-dimensional minimal lc center $C$. Suppose that $x$ is not itself an lc center. Then for a suitable dlt modification $f\colon Y\to X$, with reduced exceptional divisor $E$, one has
\[
    H_x^2(X,\mathcal O_X)
    \simeq
    H_x^0\!\left(C,R^1(f|_E)_*\mathcal O_E\right);
\]
see \cite[Corollary~18 and Remark~19]{Arvidsson-Solo}. Thus, in positive
characteristic, failure of the characteristic-zero depth bound is governed
by the occurrence of wild fibers.

We now turn to the complementary situation in which the closed point is
itself an isolated LC center. In characteristic zero, Chou
established a delicate relationship between the local cohomology, and hence
the depth, of $\mathcal O_{X,x}$ and the cohomology of the reduced
discrepancy-$-1$ divisor lying over $x$; see \cite{Chou}.

Using Theorem \ref{Thm:mainKovacs}, we
extend Chou's depth formula to log canonical threefolds under the
hypotheses of \cref{Thm:mainKovacs}. We also prove that the derived pushforward of $\cO_{B_\pi}$, and hence
its higher cohomology groups, is independent of resolution (compare \cite{Fujino-isolated,Ishii-Isolated-Gor-Sing}).

\begin{theorem}(cf. \cite{Chou})\label{thm:intro-depth}
Assume the hypotheses of \Cref{Thm:mainKovacs}, and let $x\in X$ be a closed point isolated in the non-klt locus. After localizing at $x$, let $f\colon Y\to X$ be a projective log resolution and let $E=B_f$. Thus $E$ is a reduced SNC surface and $f(E)=\{x\}$. Then
$$
 H_x^2(X,\cO_X)^\vee\simeq H^1(E,\omega_E),
$$
where $(-)^\vee$ denotes Matlis duality over $\cO_{X,x}$. Consequently,
$$
 \depth_x\cO_X\ge 3
 \quad\Longleftrightarrow\quad
 H^1(E,\cO_E)=0.
$$
\end{theorem}
\begin{theorem}(cf.\cite{Fujino-isolated,Ishii-Isolated-Gor-Sing})\label{thm:intro-intrinsic}
In the setting of \Cref{thm:intro-depth},
$H^i(E,\cO_E)$ for $i>0$
are independent of the projective log resolution.
\end{theorem}

\cite{Hacon-Alexeev} and \cite{KovacsIrr} introduced the notions of non-rational/irrational centres in char 0 and proved that for a log canonical pair $(X, B)$, every irrational centre is a non-KLT centre. In the next theorem, we provide the following counterexamples.

\begin{theorem}\label{thm:intro-examples}
The following hold.
\begin{enumerate}[label=\textnormal{(\roman*)}]
\item(cf.\cite{CT}) In characteristic $2$ there is a three-dimensional plt pair for which Kov\'acs vanishing fails.
\item(cf.\cite{Bernasconi}) In characteristic $3$ there is a normal $\Q$-factorial LC threefold pair $(X,\Delta_X)$ whose boundary has the standard coefficient $\frac12$, and a closed isolated LC center $v\in X$, such that
$
 \depth_v\cO_X=2,$ and $
 H^1(B_\pi,\cO_{B_\pi})=0
$
for the complete reduced discrepancy $-1$ divisor $F$ on a log resolution.
\item(cf.\cite{ABP,Hacon-Alexeev,KovacsIrr}) In every characteristic $p>0$, there is a LC threefold with an absolute irrational center which is not a LC center.
\end{enumerate}
\end{theorem}

\section{Acknowledgements}
 We would like to thank our advisors, Christopher Hacon, Karl Schwede, and Joe Waldron, for their guidance and encouragement. We are particularly grateful to Christopher Hacon for helpful discussions around Proposition \ref{prop:DCC-Arvidsson}, \S \ref{IndepofRes}, and \S \ref{IrrCent}. We would also like to thank Daniel Apsley, Fabio Bernasconi, Karl Schwede, and Joe Waldron for their interest in our work and valuable feedback. We are grateful to Joe Waldron for his careful reading of several drafts. The first-named author was partially supported by a Summer Research Fellowship at the University of Utah, NSF research grant DMS-2301374, and a grant from the Simons Foundation, SFI-MPS-MOV-00006719-07, and the second-named author was supported by NSF Grant DMS \#2401279 and by the Simons Foundation (award \#850684, JW) while working on this project.

Finally, we thank the organisers of the Notions of Singularity in Different Characteristics workshop at the Banff International Research Station, where the main ideas of this project were developed.
\subsection{AI Statement}
An earlier version of this preprint was run through ChatGPT Astra and Claude to find mistakes. Their main contributions were finding a gap in the hypothesis of Proposition \ref{prop:DCC-Arvidsson} (which we fixed ourselves) and pointing out \cite{ArvPosva} which we could use to simplify some of our calculations in \S  \ref{subsec:charthree}. In particular, all ideas and proofs are authors', who also take full responsibility for their correctness.
\section{Preliminaries}\label{sec:setup}
We work over an excellent local domain $R$ admitting a dualizing complex.
\subsection{Conventions and resolutions}
By a \emph{log pair} $(X,\Delta)$ we mean a normal variety $X$ and an effective $\Q$-divisor $\Delta$ with coefficients in $[0,1]$ such that $K_X+\Delta$ is $\Q$-Cartier.

Resolution of singularities for quasi-excellent schemes of dimension at
most three is due to Cossart--Piltant
\cite[Theorem~1.1]{CP19}. Projective log resolutions that we will be using here here exist by \cite[Proposition~2.14, Remark 2.15]{BMPSTWW} and
\cite[Corollary~4 and Remark~6]{KW}. A thrifty resolution is in the
sense of \cite[Definition~2.79]{KollarSingularities}.

Any two proper resolutions have a common regular model projective over each. By \cite[Corollary~4 and Remark~6]{KW}, the common model can be chosen projective over $X$. Its maps to the resolutions are then projective by \cite[Tag~0C4Q]{stacks-project}.

A coefficient is \emph{standard} if it belongs to
$$
 \left\{1-\frac1m\,\middle|\,m\in\Z_{>0}\right\}\cup\{1\}.
$$
We refer to \cite{KM,KollarSingularities} for the standard terminology.

For a DCC set $I\subseteq[\tfrac12,1]\cap\Q$, define
\[
J_I=\left\{1-\frac{1-b}{m}\ \middle|\ b\in I\cup\{0,1\},\ m\in\Z_{>0}\right\}\setminus\{0\}.
\]
For boundary coefficients in $I\cup\{1\}$, the nonzero coefficients of the different lie in $J_I$, see \cite[(3.36.1) and Corollary~3.45]{KollarSingularities}. Moreover, $J_I$ is a DCC subset of $[\tfrac12,1]$.

\subsection{Main tools}
A projective surface pair $(S,B)$ is called \emph{log Calabi--Yau} if it is log canonical and $
 K_S+B\equiv0.$

\begin{theorem}(\cite[Theorems~1.1 and~1.4]{Kaw25})\label{thm:surface-DCC}
Let
$
 I\subseteq[\tfrac12,1]\cap\Q
$
be a fixed DCC set. There exists an integer $p_0(I)>0$ with the following property. Let $k$ be a perfect field of characteristic $p>p_0(I)$, and let $(S,B)$ be a projective log Calabi--Yau surface over $k$ whose nonzero boundary coefficients belong to $I$. Let $D$ be an integral Weil divisor and let $\Theta$ be a $\Q$-divisor satisfying
$
 \Supp\Theta\subseteq\Supp B$ and $
 \floor{\Theta}=0.
$
If
$
 A:=D-(K_S+\Theta)
$
is $\Q$-Cartier, nef, and big, then
$H^i(S,\cO_S(D))=0$ for all $i>0$.
\end{theorem}
\begin{remark}
The proof is the argument of \cite[Theorem~A]{Kaw25}. The standard
coefficient hypothesis in that theorem was used only to get
log liftability. Replace this by
\cite[Theorem~1.4]{Kaw25}.
\end{remark}

We use the following special dlt modification.

\begin{proposition}(\cite[Proposition~10]{Arvidsson-Solo})\label{prop:dlt-model}
Let $(X,\Delta)$ be a quasi-projective log canonical threefold whose closed points have perfect residue fields of
characteristic greater than 5. Then there is a projective crepant birational morphism
$g\colon(Y,\Delta_Y)\longrightarrow(X,\Delta)$ such that $(Y,\Delta_Y)$ is dlt, every component of $\Exc(g)$ is $\Q$-Cartier and appears in $\Delta_Y$ with coefficient one, and there is an effective Cartier divisor $H$ with $\Supp H=\Exc(g)$ such that $-H$ is $g$-ample.
\end{proposition}
\begin{proposition}(\cite[Theorems~1.3 and~3.1(3)]{BK})\label{prop:rational-pair}
Let $(Y,\Gamma)$ be an excellent dlt threefold with a dualizing complex whose closed points have perfect residue fields of characteristic greater than $5$. If $B$ is a reduced union of components of $\floor{\Gamma}$, then $(Y,B)$ is a rational pair. In particular, for thrifty log resolution $h\colon W\to Y$, we have $Rh_*\cO_W(-h_*^{-1}B)\simeq\cO_Y(-B).$
\end{proposition}

For lc threefolds without zero-dimensional lc centers in residue characteristic $p>5$, and for standard coefficients in sufficiently large residue characteristic, \cite[Propositions~2.13 and~2.14]{ArvPosva} proves $R^1g_*\cO_Y(-E)=0$ on a suitable dlt modification. We use the restriction argument of \cite[proof of Proposition~8]{Arvidsson-Solo}.

\begin{proposition}\label{prop:DCC-Arvidsson}
Fix a DCC set
$ I\subseteq[\tfrac12,1]\cap\Q$
and let $p_I:=\max\{5,p_{0}(J_I)\}.
$ Let
$g\colon(Y,\Delta_Y)\longrightarrow(X,\Delta)$
be the dlt modification in Proposition~\ref{prop:dlt-model}. Suppose that every closed point of $X$ has perfect residue field
and that either \textup{(a)} or \textup{(b)} of
Theorem~\ref{Thm:mainKovacs} holds, with $I$ and $p_I$ as above in case \textup{(b)}. Let $D$ be an integral Weil divisor and let $0\leq\Theta\leq\Delta_Y,$ and $\floor{\Exc_g(\Theta)}=0.$
If $D-(K_Y+\Theta)$ is $\Q$-Cartier and $g$-nef, then
$
 R^ig_*\cO_Y(D)=0$ for $i>0$
\end{proposition}
\begin{proof}
We may assume that $X$ is affine and $g$ is not an isomorphism. Put $L:=D-(K_Y+\Theta)$. As in Proposition~\ref{prop:dlt-model}, choose effective integral $H$ with $\Supp H=\Exc(g)$ and $-H$ $g$-ample. Since exceptional primes are $\Q$-Cartier, we may replace $H$ by a Cartier multiple.

Fix an integer $n>0$. We religiously follow the construction in \cite[proof of Proposition~8]{Arvidsson-Solo} for $0\leq t\leq n$. Since $\lfloor\Exc_g(\Theta)\rfloor=0$, the filtration ends at $\cO_Y(D-nH)$ and every restriction occurs at $t>0$. Thus $L-tH$ is $g$-ample at each step.

For any integral exceptional $Z$, $D-Z+(\Delta_Y-\Theta)=K_Y+\Delta_Y+L-Z$ is $\Q$-Cartier. By \cite[Theorems~1.3 and~3.1(2),(4)]{BK}, $\cO_Y(D-Z)$ is Cohen--Macaulay and every exceptional primes normal.

At a step removing an exceptional prime $S$ from $D-Z$, we have $0\leq\Theta+tH-Z-S\leq\Delta_Y-S$. Apply \cite[Lemma~2.1]{BK} to $(Y,\Delta_Y)$ and $D-Z-S-K_Y$. Since $\cO_Y(D-Z-S)$ is $S_3$ along $S$, we get
\[
0\longrightarrow\cO_Y(D-Z-S)\longrightarrow\cO_Y(D-Z)\longrightarrow\cO_S(D_S)\longrightarrow0
\]
with $D_S$ integral. Since $D-Z-S-K_Y=(L-tH)+(\Theta+tH-Z-S)$, the residue formula
gives $D_S\sim_\Q K_S+(L-tH)|_S+\Theta_S$. At each codimension-one point of $S$, the surface is $\Q$-factorial by \cite[Lemma~2.30]{ABP}. Choose something linearly equivalent to $D-Z-S-K_Y$ whose support does not contain $S$. By \cite[Lemma~2.32(i) and the proof of Proposition~2.33]{ABP}, the correction to its restriction lies in between zero and Diff$_S(0)$. Thus we may choose $\Theta_S$ so that $0\leq\Theta_S-(\Theta+tH-Z-S)|_S\leq$Diff$_S(0).$

Since $0\leq\Theta+tH-Z-S\leq\Delta_Y-S$, it follows that $0\leq\Theta_S\leq(\Delta_Y-S)|_S+$Diff$_S(0)
=$Diff$_S(\Delta_Y-S)$.

$(S,$ Diff$_S(\Delta_Y-S))$ is dlt by adjunction. So $S$ is $\Q$-factorial by \cite[Corollary~4.11]{Tanaka}. As in \cite[Remark~6]{Arvidsson-Solo}, replace $\Theta_S$ by $(1-\varepsilon)\Theta_S$ for $0<\varepsilon\ll1$. Then $(S,(1-\varepsilon)\Theta_S)$ is klt and $D_S-(K_S+(1-\varepsilon)\Theta_S)$ is ample over $g(S)$ \cite[Remark~6]{Arvidsson-Solo}.

If $g(S)$ is a curve, then $R^i(g|_S)_*\cO_S(D_S)=0$ for $i>0$ by \cite[Lemma~12]{Arvidsson-Solo}, after normalizing.

If $g(S)=x$, then $x$ is an lc center. So we are in case (b). By crepantness and adjunction, $(S,$Diff$_S(\Delta_Y-S))$ is a projective log Calabi--Yau surface over the perfect field $\kappa(x)$, with coefficients in $J_I$. Theorem~\ref{thm:surface-DCC} then gives $H^i(S,\cO_S(D_S))=0$ for $i>0$.

Now, by Serre, $R^ig_*\cO_Y(D-nH)=0$ for $i>0$ and $n\gg0$. Reverse the restriction steps and get $R^ig_*\cO_Y(D)=0$ for $i>0$.
\end{proof}

\section{Kov\'acs vanishing and isolated centers}\label{sec:main}
Take
$ g\colon(Y,\Delta_Y)\longrightarrow(X,\Delta)$
the crepant dlt modification coming from Proposition~\ref{prop:dlt-model}. Let
$
 B:=\floor{\Delta_Y},$ and $
 \Theta:=\Delta_Y-B.$ By \cite[Theorem~2.13 and Remark~2.15]{BMPSTWW}, we may choose a proper log resolution $h\colon W\to (Y,\Delta_Y)$ which is an isomorphism over the SNC locus of $(Y,\Delta_Y)$. Since $(Y,\Delta_Y)$ is dlt, this resolution is thrifty for $(Y,B)$, \cite[Definition~2.79]{KollarSingularities}. Let
\[
    \pi_0:=g\circ h,\quad \mathrm{and} \quad F:=h_*^{-1}B.
\]
\begin{lemma}\label{lem:boundary-identification}
We have $F=B_{\pi_0}$ and $R^i\pi_{0*}\cO_W(-F)=0$ for $i>0$.
\end{lemma}
\begin{proof}
The dlt pair $(Y,\Delta_Y)$ is SNC at the generic point of every lc center. By construction, $h$ is an isomorphism there. Hence no $h$-exceptional divisor has discrepancy $-1$. Since $g$ is crepant, $B_{\pi_0}=h_*^{-1}\floor{\Delta_Y}=F$. So we get the first claim.

By Proposition~\ref{prop:rational-pair},
$
 Rh_*\cO_W(-F)\simeq\cO_Y(-B).
$
So, applying $Rg_*$, we get $R\pi_{0*}\cO_W(-F)\simeq Rg_*\cO_Y(-B)$. Note $-B-(K_Y+\Theta)=-(K_Y+\Delta_Y)=-g^*(K_X+\Delta)$ which is $\Q$-Cartier and $g$-numerically trivial. Also $0\le\Theta\le\Delta_Y$ with $\Theta$ having no $g$-exceptional components. So $R^ig_*\cO_Y(-B)=0$ for $i>0$ by Proposition~\ref{prop:DCC-Arvidsson}. Hence we get the second claim.
\end{proof}

\subsection{Independence of resolutions}\label{IndepofRes}
Let $q\colon V\longrightarrow X$ be a proper log resolution. Take
$K_V+\Gamma_V=q^*(K_X+\Delta)$, and $B_V:=\Gamma_V^{=1}.$

\begin{lemma}\label{lem:coefficient-one-reduction}
Let $\sigma\colon W\to V$ be projective and birational, with $W$ regular. Then the crepant pullbacks of $(V,\Gamma_V)$ and $(V,B_V)$ to $W$ have the same components of coefficient one.
\end{lemma}
\begin{proof}
Let $E$ be a divisor over $V$ with discrepancy $-1$ for either pair.
If its center were contained in a component of $\Gamma_V-B_V$, that
component would have coefficient less than one in both boundaries.
Replacing the boundary by the reduced divisor $\Supp\Gamma_V$ would
then lower the discrepancy of $E$ below $-1$. This is impossible
because $(V,\Supp\Gamma_V)$ is a SNC pair and
hence log canonical. Thus $\Gamma_V$ and $B_V$ agree at the generic
point of the center of $E$. So $E$ has discrepancy $-1$ for both pairs.
\end{proof}

The next lemma should be well-known to experts; compare
\cite[Remark~4]{Fujino-Kovacsvanishing} and \cite[Section~4]{Chou}.
\begin{lemma}\label{lem:one-blowup}
Let
$
 \beta\colon V^+\longrightarrow V$
be the blow-up of a regular irreducible center $C$ of codimension $c\geq2$ having simple normal crossings with $\Supp\Gamma_V$. Then
$
 R\beta_*\cO_{V^+}(-B_{V^+})\simeq\cO_V(-B_V).
$
\end{lemma}
\begin{proof}
Let \(F\) be the exceptional divisor, and \(D_1,\ldots,D_s\) be the components of \(\Gamma_V\) containing the generic point of \(C\), with coefficients \(b_1,\ldots,b_s\). Let $r$ be the number of components of $B_V$ containing $C$. The SNCness gives $s\le c$. The coefficient of $F$ in $\Gamma_{V^+}$ is $\sum_{j=1}^s b_j-(c-1)$. Since each $b_j\le1$, this coefficient is one exactly when $r=c$. The divisors $\beta^*B_V$ and $B_{V^+}$ agree away from $F$. Their difference has coefficient $r$ along $F$ if $r<c$, and $c-1$ if $r=c$. Thus $\beta^*B_V-B_{V^+}=mF$ with $m=\min\{r,c-1\}$.

The map $F\to C$ is a $\mathbb P^{c-1}$-bundle with $\cO_F(F)=\cO_F(-1)$. For $1\le m\le c-1$, all direct images of $\cO_F(-m)$ vanish \cite[Tag~01XX]{stacks-project}. From
\[
0\to\cO_{V^+}((m-1)F)
\to\cO_{V^+}(mF)
\to\cO_F(-m)\to0
\]
and \(R\beta_*\cO_{V^+}\simeq\cO_V\) (\cite[Theorem~1.1]{CR}), we get $R\beta_*\cO_{V^+}(mF)\simeq\cO_V$ for $0\le m\le c-1$ by induction.
Projection formula finishes the claim.
\end{proof}
\begin{lemma}\label{lem:zero-pushforward}
Let $\sigma\colon W\to V$ be a projective birational morphism between regular log resolutions over $X$. Then
$
 \sigma_*\cO_W(-B_W)=\cO_V(-B_V).
$
\end{lemma}
\begin{proof}
By Lemma~\ref{lem:coefficient-one-reduction}, $B_W$ is the sum of the prime divisors on $W$ with discrepancy $-1$ for $(V,B_V)$. It contains the strict transform of $B_V$. Under $\sigma_*\cO_W=\cO_V$, this gives $\sigma_*\cO_W(-B_W)\subseteq\cO_V(-B_V)$.

$(V,B_V)$ is klt off $B_V$. So every component of $B_W$ maps into $B_V$. Since $B_V$ is Cartier and $B_W$ is reduced, $B_W\le \sigma^*B_V$. By the projection formula,
\[
\cO_V(-B_V)=\sigma_*\cO_W(-\sigma^*B_V)\subseteq\sigma_*\cO_W(-B_W).
\]
\end{proof}
Now we are ready to prove Kov\'acs vanishing!
\begin{proof}[Proof of Theorem~\ref{Thm:mainKovacs}]
Let $q\colon V\to X$ be a proper log resolution. Take a projective birational morphism $\sigma\colon W\to V$ such that $q\sigma$ is a log resolution. Apply \cite[2.10 and proof of Theorem~2.12]{BK} to $\sigma$ with $\Exc(q)\cup\Supp(q_*^{-1}\Delta)$ as boundary on $V$. This gives a morphism $\rho\colon Z\to W$ such that $\tau=\sigma\rho$ is a composition of admissible blow-ups. At each step, the center has SNC with the union of the exceptional divisor over $X$ and the strict transform of $\Supp\Delta$. So $q\tau$ is a log resolution. Since $\tau$ is projective and $\sigma$ separated, $\rho$ is projective \cite[Tag~0C4Q]{stacks-project}. Apply Lemma~\ref{lem:one-blowup} at each step. We get $R\tau_*\cO_Z(-B_Z)\simeq\cO_V(-B_V)$.

We prove $R^i\sigma_*\cO_W(-B_W)=0$ for all such $\sigma$ by induction on $i>0$. By induction, $R^j\rho_*\cO_Z(-B_Z)=0$ for $0<j<i$. $\rho_*\cO_Z(-B_Z)=\cO_W(-B_W)$ by Lemma~\ref{lem:zero-pushforward}. So the edge map therefore gives $R^i\sigma_*\cO_W(-B_W)\hookrightarrow R^i\tau_*\cO_Z(-B_Z)=0$. Lemma~\ref{lem:zero-pushforward} also gives $\sigma_*\cO_W(-B_W)=\cO_V(-B_V)$. So,
\begin{equation}\label{eq:resolution-comparison}
\cO_V(-B_V)\xrightarrow{\sim}R\sigma_*\cO_W(-B_W).
\end{equation}
These maps are compatible with composition because in degree zero they are pullbacks of functions.

Let $\pi\colon Y\to X$ be a proper log resolution. Take a log resolution that dominates both $\pi$ and $\pi_0$ by projective morphisms. By \eqref{eq:resolution-comparison}, all three resolutions give the same higher direct images. They vanish by Lemma~\ref{lem:boundary-identification}.

Each component of $B_\pi$ has discrepancy $-1$. Over $X\setminus\pi(B_\pi)$, the crepant pair on $Y$ is sub-klt. Thus $\pi(B_\pi)$ is the non-klt locus. By normality, $\pi_*\cO_Y=\cO_X$. Since $B_\pi$ is reduced, we get $R\pi_*\cO_Y(-B_\pi)$ is canonically the ideal sheaf of the reduced non-klt locus.
\end{proof}

\subsection{Isolated log canonical centers}
First we recall the following weak Grauert--Riemenschneider vanishing (\cite{Arvidsson-Solo}).

\begin{corollary}\label{cor:weak-GR}
Suppose that every closed point of $X$ has perfect residue field and that either condition \textnormal{(a)} or \textnormal{(b)} of \Cref{Thm:mainKovacs} holds. Then every proper log resolution $f\colon Y\to X$ satisfies
\[
 R^if_*\omega_Y=0\quad \text{for all} \ i>0.
\]
\end{corollary}
\begin{proof}
On the special dlt model $g\colon(V,\Gamma)\to(X,\Delta)$, apply Proposition~\ref{prop:DCC-Arvidsson} with $D=K_V$ and $\Theta=0$. So $R^ig_*\omega_V=0$ for $i>0$.

Take the good log resolution $h\colon W\to V$ as before. $V$ is Cohen--Macaulay by
\cite[Theorems~1.3 and~3.1(1)]{BK}.
Proposition~\ref{prop:rational-pair} with $B=0$ gives $Rh_*\omega_W\simeq\omega_V$, by duality. So the assertion holds for $gh$. Compare $gh$ and $f$ on a common regular model projective over both. For each upper map $u$, \cite[Theorem~1.1]{CR} gives $Ru_*\omega=\omega$. The claim follows by composition.
\end{proof}
\begin{proof}[Proofs of Theorems~\ref{thm:intro-depth} and~\ref{thm:intro-intrinsic}]
For the first isomorphism, the proof is Chou's argument
\cite[Theorem~3.1 and Corollary~3.2]{Chou}. Theorem~\ref{Thm:mainKovacs} and
Corollary~\ref{cor:weak-GR} provide, respectively, the Kov\'acs and
Grauert--Riemenschneider vanishing inputs in our setting. $X$ is Cohen--Macaulay at $x$ if and only if
$
 H^1(E,\cO_E)=0
$ since the local ring has dimension three. Finally, we show $
 Rf_*\cO_E
$
is independent of the projective log resolution, and so is
$ H^i(E,\cO_E)$ for $i>0
$. Let \(\sigma\colon W\to Y\) compare two resolutions. The inclusion $
\sigma^*\cO_Y(-B_Y)\longrightarrow\cO_W(-B_W)$
gives a morphism from $$
0\to\cO_Y(-B_Y)\to\cO_Y\to\cO_{B_Y}\to0$$
to the derived pushforward of
\[
0\to\cO_W(-B_W)\to\cO_W\to\cO_{B_W}\to0.
\]
The first two vertical maps are isomorphisms by \eqref{eq:resolution-comparison} and \cite[Theorem~1.1]{CR}. Hence $
\cO_{B_Y}\xrightarrow{\sim}R\sigma_*\cO_{B_W}.$
These isomorphisms are compatible with composition. Passing to a common resolution therefore shows that $Rf_*\cO_{B_f}$ is independent of the resolution. If \(B_f\) lies over \(x\), then
$
(R^if_*\cO_{B_f})_x=H^i(B_f,\cO_{B_f}),$
which gives the claim.
\end{proof}

\section{Counterexamples}\label{sec:examples}
\subsection{Characteristic two}
Now we use Cascini--Tanaka's non-normal plt center example \cite[Theorem~1.1]{CT} to show that Kov\'acs vanishing might fail in low char.
\begin{proposition}\label{prop:char2}
In characteristic two there is a $\Q$-factorial plt threefold $(Z,E)$ such that, for every projective log resolution
$h\colon W\to Z$ and $F=B_h$, we have $R^1h_*\cO_W(-F)\ne0.$
\end{proposition}
\begin{proof}
Let $(Z,E)$ be the $\Q$-factorial plt pair of \cite[Theorem~1.1]{CT} with the prime $E$ not normal. Choose a log resolution $h\colon W\to Z$, and let $F$ be the strict transform of $E$. Since the pair is plt, $F=B_h$.

Push forward
\begin{equation}\label{eq:char2-sequence}
 0\longrightarrow\cO_W(-F)
 \longrightarrow\cO_W
 \longrightarrow\cO_F
 \longrightarrow0.
\end{equation}
The variety $Z$ is normal. So $h_*\cO_W=\cO_Z$. If
$
 \nu\colon E^\nu\longrightarrow E
$
is the normalization, then $h|_F$ factors through $\nu$ and
$ h_*\cO_F=\nu_*\cO_{E^\nu}.
$
The map $\cO_Z\to h_*\cO_F$ factors through
$\cO_E\longrightarrow\nu_*\cO_{E^\nu}.
$
This map is not surjective because $E$ is not normal. Its nonzero cokernel injects, by the long exact sequence associated to \eqref{eq:char2-sequence}, into
$
 R^1h_*\cO_W(-F).
$
\end{proof}

\subsection{Characteristic three}\label{subsec:charthree}
Let $k$ be an algebraically closed field of characteristic three. We use Bernasconi's surface from \cite[Section~3.1]{Bernasconi}. Let $T$ and $E_1,E_2,E_3$ be the surface and prime divisors of \cite[Section~3.1]{Bernasconi}. By \cite[Proposition~3.2]{Bernasconi}, $T$ is a klt del Pezzo surface with $\rho(T)=1$. It is $\Q$-factorial by \cite[Corollary~4.11]{Tanaka}.

Let $A=E_2+E_3-E_1$. By \cite[Theorem~3.6]{Bernasconi}, $A$
is ample and
\begin{equation}\label{eq:Bernasconi-nonvanishing}
 H^1(T,\cO_T(-A))\neq0.
\end{equation}
Put $\Delta_T=\tfrac12(E_2+E_3)$.
\begin{lemma}\label{lem:surface-klt}
The pair $(T,\Delta_T)$ is klt and $K_T+\Delta_T\sim_\Q0$.
\end{lemma}
\begin{proof}
By \cite[Lemma~4.4]{ArvPosva}, $E_1\sim_\Q-K_T$ and $(T,E_1)$ is plt. The three points in \cite[Section~3.1]{Bernasconi} can be numbered in any order. So the same holds for $E_2$ and $E_3$. In particular, $\Delta_T\sim_\Q-K_T$. For every exceptional prime divisor $F$ over $T$,
\[
a(F,T,\Delta_T)=\frac{a(F,T,E_2)+a(F,T,E_3)}{2}>-1.
\]
Since $\floor{\Delta_T}=0$, the pair $(T,\Delta_T)$ is klt.
\end{proof}

Consider the cone
$X_0=\Spec\bigoplus_{m\geq0}H^0(T,\cO_T(mA))$with vertex $v$ and birational model
\[
 Y_0=\Spec_T\!\left(\bigoplus_{m\geq0}\cO_T(mA)\right),
 \quad f_0\colon Y_0\longrightarrow X_0,\quad
 \pi\colon Y_0\longrightarrow T .
\]
The exceptional divisor $E$ of $f_0$ is the zero section and $E\simeq T$. Moreover $Y_0\setminus E\simeq X_0\setminus\{v\}$ \cite[Section~2.3]{Bernasconi}. Let $\pi^*E_i$ be the reduced closure in $Y_0$ of $\pi^{-1}(E_i\cap T_{\rm reg})$. Take $\Delta_{X_0}=f_{0*}\pi^*\Delta_T$ for $\pi^*\Delta_T=(\pi^*E_2+\pi^*E_3)/2$ .

\begin{lemma}\label{lem:cone-properties}
 $X_0$ is normal and $\Q$-factorial.  $\Delta_{X_0}$ has two prime components each with coefficient $\frac12$. With compatible canonicals, we have
\[
 K_{Y_0}+E+\pi^*\Delta_T=f_0^*(K_{X_0}+\Delta_{X_0}).
\]
Finally $(Y_0,E+\pi^*\Delta_T)$ is plt near $E$.
\end{lemma}
\begin{proof}
The cone $X_0$ is normal $\Q$-factorial by \cite[Theorem~4.2]{Bernasconi}. $\pi^*E_2$ and $\pi^*E_3$ have distinct images because $f_0$ is an isomorphism off $E$.

By Lemma~\ref{lem:surface-klt}, $K_T+\Delta_T\sim_\Q 0$. The crepant equality follows from \cite[(2.4)]{ArvPosva}. Near every codimension-one point of $T$, the divisor $A$ is Cartier and $\pi$ is the projection of a line bundle. Adjunction at these points gives Diff$_E(\pi^*\Delta_T)=\Delta_T$ under $E\simeq T$. This determines the different as a divisor on $E$. Since $(T,\Delta_T)$ is klt, inversion of adjunction shows that $(Y_0,E+\pi^*\Delta_T)$ is plt near $E$ \cite[Lemma~2.26]{ArvPosva}.
\end{proof}
\begin{theorem}\label{thm:char3}
There is a normal $\Q$-factorial three-dimensional log canonical pair $(X,\Delta_X)$ over $k$ whose only nonzero boundary coefficient is $\frac12$, and a closed isolated LC center $v\in X$, such that
$\depth_v\cO_X=2.$ If $h\colon W\to X$ is any projective log resolution of the germ at $v$ and $F = B_h$, then $F$ is a smooth projective rational surface. In particular, $H^1(F,\cO_F)=0$.
\end{theorem}
\begin{proof}
Choose an open neighborhood \(U\subset Y_0\) of \(E\) on which \((Y_0,E+\pi^*\Delta_T)\) is plt. Note that \(f_0\) is proper. Choose \(r>0\) such that \(rA\) is Cartier and the \(r\)-th Veronese is generated in degree one. Then \(Y_0\) and \(X_0\) are finite over their \(r\)-th Veroneses, and $Y_0^{(r)}\longrightarrow X_0^{(r)}$ is projective. Hence \(Y_0\to X_0^{(r)}\) is proper. Since \(X_0\to X_0^{(r)}\) is finite and separated, \(f_0\colon Y_0\to X_0\) is proper.

Let $X=X_0\setminus f_0(Y_0\setminus U)$. This is an open neighborhood of $v$ because $f_0$ is proper and $(f_0^{-1}(v))_{\rm red}=E$. Let $\Delta_X:=\Delta_{X_0}|_X.$ The crepant equality
\[
\left.(K_{Y_0}+E+\pi^*\Delta_T)\right|_{f_0^{-1}(X)}
=
\bigl(f_0\vert_{f_0^{-1}(X)}\bigr)^*(K_X+\Delta_X).
\]
shows that \((X,\Delta_X)\) is lc and with \(E\) having discrepancy \(-1\). Thus \(v\) is an lc center. Since the pair upstairs is plt and \(f_0\) is an isomorphism away from \(E\), \((X\setminus\{v\},\Delta_X)\) is klt. Hence \(v\) is an isolated lc center.

By \cite[Proposition~2.6]{Bernasconi},
\[
H_v^2(X_0,\cO_{X_0})
\simeq
\bigoplus_{m\in\Z}H^1(T,\cO_T(mA)).
\]
Degree $(-1)$ piece is nonzero by \eqref{eq:Bernasconi-nonvanishing}. Therefore
$
H_v^2(X,\cO_X)\neq0.
$
So \(\depth_v\cO_X\le2\). Since \(X\) is normal, it is \(S_2\). Hence \(\depth_v\cO_X\ge2\). Thus $\depth_v\cO_X=2.$

Let $h\colon W\to X$ be any projective log resolution. Since the germ is lc, the crepant boundary on $W$ has a component of coefficient one. The plt model shows that the valuation of $E$ is the only lc place. Hence there is a unique component $F$ of coefficient one. Since $F$ and $E$ define the same valuation, they are birational.

Now \(F\) is a smooth projective surface birational to \(E\simeq T\). By \cite[Lemma~5.1]{Bernasconi}, \(T\) is rational. Hence \(F\) is rational, and
$
H^1(F,\cO_F)=0.$
\end{proof}

This example shows that the characteristic hypothesis in Theorem~\ref{thm:intro-depth} cannot be omitted, even for a $\Q$-factorial pair with a standard coefficient.

\subsection{Irrational centers}\label{IrrCent}
Following \cite{Hacon-Alexeev,KovacsIrr}, a closed irreducible subvariety $Z\subseteq X$ is an \emph{absolute irrational center} if, for every proper resolution $f\colon Y\to X$, the generic point of $Z$ is an associated point of some $R^if_*\cO_Y$, $i>0$.

\begin{proposition}\label{prop:irrational}
In every positive characteristic there is a log canonical threefold with an absolute irrational centre which is not a log canonical centre.
\end{proposition}
\begin{proof}
Fix a prime $p$ and let $C$ be the spectrum of an equicharacteristic DVR with algebraically closed residue field of characteristic $p$. Choose a regular minimal elliptic fibration
$
 q\colon S\longrightarrow C
$
with a wild fiber as in \cite[Example~4.7]{ABP}. Write
$$
 R^1q_*\cO_S=M\oplus T,
$$
where $M$ is locally free and $0\ne T$ is torsion. Let $x$ be a closed point in $\Supp T$. Choose a relatively ample invertible sheaf $L$ on $S$.

Let
\[
 X=\Spec_C\!\left(\bigoplus_{m\geq0}q_*L^{\otimes m}\right),
 \qquad
 Y=\Spec_S\!\left(\bigoplus_{m\geq0}L^{\otimes m}\right),
\]
and let $f\colon Y\to X$ be the natural morphism. By \cite[Proposition~5.3]{ABP}, $X$ is a log canonical threefold, $Y$ is regular, and the relative vertex
$
 V_C\simeq C
$
is the unique LC center. In particular, the closed point $x\in V_C$ is not an lc center.

By \cite[Claim~5.6(ii)]{ABP},
$
 R^1f_*\cO_Y\simeq i_*(M\oplus T),
$
where $i\colon V_C\hookrightarrow X$ is the inclusion. Since $T_x$ has nonzero finite length, there is a nonzero element annihilated by the maximal ideal of $\cO_{C,x}$. Its image in
$(i_*T)_x$ has annihilator $\mathfrak m_{X,x}$. Thus $x$ is an associated point of $R^1f_*\cO_Y$.

Let $f'\colon Y'\to X$ be any proper resolution. Choose a common regular model $W$ projective over $Y$ and $Y'$. By \cite[Theorem~1.1]{CR}, the derived pushforward of $\cO_W$ along either upper map is the structure sheaf of its target. Hence $Rf'_*\cO_{Y'}\simeq Rf_*\cO_Y$. It follows that $x$ is an absolute irrational center.
\end{proof}
\begin{proof}[Proof of Theorem~\ref{thm:intro-examples}]
Part (i) is Proposition~\ref{prop:char2}, part (ii) is Theorem~\ref{thm:char3}, and part (iii) is Proposition~\ref{prop:irrational}.
\end{proof}

\bibliographystyle{alpha}
\bibliography{main}
\end{document}